\documentclass[11pt,openbib]{article} 
\usepackage{amsfonts,amsthm}
\usepackage{amssymb,amsmath}

\newcommand{\kk}{\mathrm{k}}

\newcommand{\Z}{\mathbb{Z}}

\newcommand{\imm}{\mathrm{im}\, }

\newcommand{\setrule}{\, \mathop{\rule[-4pt]{.5pt}{13pt}\, }\nolimits}

\newcommand{\spann}{\mathop{\rm span}\nolimits}

\allowdisplaybreaks

\begin{document}
\begin{center}
{\LARGE \bf The uniform normal form} \\
\mbox{}\vspace{-.1in}  \\
{\LARGE \bf  of a linear mapping}
\end{center}
\begin{center}
\mbox{}\vspace{-.1in} \\ 
Richard Cushman\footnotemark 
\mbox{}\vspace{.35in} \\
Dedicated to my friend and colleague Arjeh Cohen on his retirement
\end{center}
\footnotetext{Department of Mathematics and Statistics, 
University of Calgary, \\ 
Calgary, Alberta, Canada, T2N 1N4.}  
\vspace{.25in}

Let $V$ be a finite dimensional vector space over a field 
$\kk $ of characteristic $0$. Let $A:V \rightarrow V$ be a linear 
mapping of $V$ into itself with characteristic polynomial 
${\chi }_A$. The goal of this paper is to give a normal 
form for $A$, which yields a better description of its structure than the classical companion matrix. This normal form does not use a factorization of ${\chi }_A$ and requires only 
operations in the field $\kk $ to compute.

\section{Semisimple linear mappings}
\label{sec1}

We begin by giving a well known criterion to determine if the 
linear mapping $A$ is semisimple, that is, every $A$-invariant subspace of $V$ has 
an $A$-invariant complementary subspace. \medskip 

Suppose that we can factor ${\chi }_A$, that is, find monic irreducible 
polynomials ${ \{ {\pi }_i \} }^m_{i=1}$, which are pairwise 
relatively prime, such that ${\chi }_A = {\prod }^m_{i=1} {\pi }^{n_i}_i$, 
where $n_i \in {\Z}_{\ge 1}$. Then 
\begin{displaymath}
{\chi}'_A =\sum^m_{j=1} (n_j {\pi }^{n_j-1}_j {\pi }'_j) 
{\prod}_{i\ne j}{\pi }^{n_i}_i   =
\big( {\prod }^m_{\ell =1} {\pi }^{n_{\ell }-1}_{\ell } \big) \big( 
\sum^m_{j=1} (n_j {\pi }'_j) \prod_{i\ne j}{\pi }_i \big) . 
\end{displaymath}
Therefore the greatest common divisor ${\chi }_A$ and its 
derivative ${\chi }'_A$ is the polynomial 
$d= {\prod }^m_{\ell =1} {\pi }^{n_{\ell }-1}_{\ell }$. The polynomial $d$ can be 
computed using the \linebreak 
Euclidean algorithm. Thus the square free factorization 
of ${\chi }_A$ is the polynomial $p = 
{\prod }^m_{\ell =1} {\pi }_{\ell } = {\chi }_A/d $, which can be 
computed without knowing a factorization of ${\chi }_A$.  \medskip 

The goal of the next discussion is to prove \medskip 

\noindent \textbf{Claim 1.1} The linear mapping $A:V \rightarrow V$ is semisimple if $p(A) =0 $ on $V$. \medskip 

Let $p = {\prod }^m_{j =1} {\pi }_j$ be the square free factorization 
of the characteristic polynomial ${\chi }_A $ of $A$. We now decompose $V$ into $A$-invariant subspaces. For each $1\le j \le m$ let $V_j = 
\{ v \in V \setrule \, {\pi }_j (A)v = 0 \} $. Then $V_j$ is an 
$A$-invariant subspace of $V$. For if $v \in V_j$, then 
${\pi }_j(A)Av = A{\pi }_j(A)v =0$, that is, $Av \in V_j$. The 
following argument shows that $V = {\bigoplus }^m_{j=1} V_j$. 
Because for $1 \le j \le m$ the polynomials ${\prod }_{i\ne j}{\pi}_i$ 
are pairwise relatively prime, there are polynomials $f_j$, $1 \le j \le m$ 
such that $1 = \sum^m_{j=1} f_j \, \big( \prod_{i\ne j}{\pi }_i \big)$. 
Therefore every vector $v \in V$ can be written as 
\begin{displaymath}
v = \sum^m_{j=1} f_j(A) \, \big( \prod_{i\ne j}{\pi }_i(A)v \big) = 
\sum^m_{j=1}f_j(A)v_j. 
\end{displaymath}
Since ${\pi }_j(A) \big( \prod_{i\ne j}{\pi }_i(A)v \big) = p(A)v =0$, 
the vector $v_j \in V_j$. Therefore $V = \sum^m_{j=1}V_j$. If for 
$i \ne j$ we have $w \in V_i \cap V_j$, then for some polynomials 
$F_i$ and $G_j$ we have $ 1 = F_i {\pi }_i + G_j {\pi }_j$, because 
${\pi }_i$ and ${\pi }_j$ are relatively prime. Consequently, 
$w = F_i(A){\pi }_i(A)w +G_j(A){\pi }_j(A)w =0$. So $V = \sum^m_{j=1} \oplus V_j$. \hfill $\square $ \medskip 

We now prove \medskip 

\noindent \textbf{Lemma 1.2} For each $1 \le j \le m$ there is a 
basis of the $A$-invariant subspace $V_j$ such that that matrix of 
$A$ is block diagonal. \medskip 

\noindent \textbf{Proof.} Let $W$ be a minimal dimensional proper 
$A$-invariant subspace of $V_j$ and let $w$ be a nonzero vector in $W$. 
Then there is a minimal positive integer $r$ such that 
$A^r w \in {\spann }_{\kk }\{ w, \, Aw, \, \ldots , \, A^{r-1}w \} = U$. 
We assert: the vectors ${\{ A^iw \} }^{r-1}_{i=0}$ are linearly 
independent. Suppose that there are $a_i \in \kk $ for $1 \le i \le r-1$ such that 
$0 = a_0w + a_1Aw + \cdots + a_{r-1}A^{r-1}w$. 
Let $t \le r-1$ be the largest index such that $a_t \ne 0$. So 
$A^tw = -\frac{a_{t-1}}{a_t} A^{t-1}w - \cdots - \frac{a_0}{a_t} w$, 
that is, $A^t w \in {\spann }_k \{ w, \, \ldots , \, A^{t-1}w \} $ and 
$t <r $. This contradicts the definition of the integer $r$. Thus the 
index $t$ does not exist. Hence $a_i =0$ for every  $0 \le i \le r-1$, 
that is, the vectors ${\{ A^iw \} }^{r-1}_{i=0}$ are linearly 
independent. 

The subspace $U$ of $W$ is $A$-invariant, for 
\begin{align}
A(\sum^{r-1}_{j=0} b_j A^jw) & = \sum^{r-2}_{j=0} b_j A^{j+1}w 
+b_{r-1}A^r w, \, \, \, \mbox{where $b_j \in \kk$} \notag \\
&\hspace{-.5in} = \sum^{r-1}_{j=1}b_{j-1}A^jw +b_{r-1}(\sum^{r-1}_{\ell =0} 
a_{\ell}A^{\ell }w), \quad \mbox{since $A^rw \in U$} \notag \\
& \hspace{-.5in} = b_{r-1}a_0w +\sum^{r-1}_{j=1}(b_{j-1}+b_{r-1}a_j)A^jw \in U. \notag 
\end{align}    
 
Next we show that there is a monic polynomial $\mu $ of degree $r$ 
such that $\mu (A) =0$ on $U$. With respect the basis 
${\{ A^iw \} }^{r-1}_{i=0}$ of $U$ we can write $A^rw = 
-a_0w - \cdots - a_{r-1}A^{r-1}w$. So $\mu (A)w =0$, where 
\begin{equation}
\mu (\lambda ) = a_0 +a_1\lambda + \cdots + a_{r-1}{\lambda }^{r-1} 
+{\lambda }^r .
\label{eq-twos1}
\end{equation}
Since $\mu (A)A^iw = A^i(\mu (A)w) =0$ for every $0 \le i \le r-1$, it 
follows that $\mu (A) = 0$ on $U$. 

By the minimality of the dimension 
of $W$ the subspace $U$ cannot be proper. But $U \ne \{ 0 \}$, since 
$w \in U$. Therefore $U = W$. Since $U \subseteq V_j$, we obtain 
${\pi }_j(A) u =0 $ for every $u \in U$. Because ${\pi }_j$ is irreducible, 
the preceding statement shows that ${\pi }_j$ is the minimum polynomial 
of $A$ on $U$. Thus ${\pi }_j$ divides $\mu $. Suppose that 
$\deg {\pi }_j =s < \deg \mu = r$. Then $A^s u^{\prime } \in 
{\spann }_k \{ u^{\prime }, \, \ldots \, A^{s-1}u^{\prime } \} = Y$ 
for some nonzero vector $u^{\prime }$ in $U$. By minimality, $Y = U$. 
But $\dim Y = s < \dim U = r$, which is a contradiction. Thus 
${\pi }_j = \mu $. Note that the matrix of $A|U$ with respect to the 
basis ${\{ A^iw \} }^{r-1}_{i=0}$ is the $r \times r$ companion matrix 
\begin{equation}
C_r = \mbox{\footnotesize $ \left( \begin{array}{lllcc}
0    & \cdots       & \cdots       & 0  &  -a_0           \\
1    & 0      &  \cdots       & 0  & -a_1           \\
\vdots    & 1 & \ddots &  \vdots   & \vdots          \\
\vdots      &           & \ddots &    0  &   -a_{r-2}  \\
0 & \cdots    &  \cdots    &  1 &  -a_{r-1}   \end{array} \right) $,}
\label{eq-twostars1}
\end{equation}
where ${\pi }_j = a_0 +a_1\lambda + \cdots + a_{r-1}{\lambda }^{r-1} 
+{\lambda }^r$. 

Suppose that $U \ne V_j$. Then there is a nonzero vector $w' \in 
V_j \setminus U$. Let $r'$ be the smallest positive integer such that 
$A^{r'}w' \in {\spann }_{\kk } \{ w', \, Aw', \, \ldots , \,$ 
$A^{r'-1}w' \} = U'$. Then by the argument in the preceding paragraph, $U'$ is a minimal $A$-invariant subspace of $V_j$ of dimension $r' =r$, whose minimal polynomial is ${\pi }_j$. Suppose that $U' \cap U \ne \{ 0 \} $. Then $U' \cap U$ is a proper $A$-invariant subspace of $U'$. By minimality $U' \cap U = U'$, 
that is, $U \subseteq U'$. But $r = \dim U = \dim U' = r'$. So $U = U'$. 
Thus $w'\in U'$ and $w' \notin U$, which is a contradiction. Therefore 
$U' \cap U = \{ 0 \} $. If $U \oplus U' \ne V_j$, we repeat the above 
argument. Using $U \oplus U'$ instead of $U$, after a finite number of repetitions we have 
$V_j = \sum^{\ell }_{i=1}\oplus U_i$, where for every $0 \le i \le \ell $ the subspace $U_i$ of 
$V_j$ is  $A$-invariant  with basis ${\{ A^k u_i \} }^{r-1}_{k =0}$ and the minimal polynomial 
of $A|U_i$ is ${\pi }_j$. With respect to the basis 
${\{ A^k u_i \}}^{(\ell , r-1)}_{(i,k) =(1,0)}$ of $V_j$ the matrix of 
$A$ is $\mathrm{diag}(C_r, \ldots , C_r)$, 
which is block diagonal. \hfill $\square $ \medskip 

For each $1 \le j \le m $ applying lemma 1.2 to $V_j$ and using 
the fact that $V = \sum^m_{j=1} \oplus V_j$ we obtain \medskip 

\noindent \textbf{Corollary 1.3} There is a basis of $V$ such that 
the matrix of $A$ is block diagonal. \medskip 

\noindent \textbf{Proof of claim 1.1} Suppose that $U$ is an $A$-invariant 
subspace $V$. Then by corollary 1.3, there is a basis ${\varepsilon }_U$ 
of $U$ such that the matrix of $A|U$ is block diagonal. By corollary 
1.3 there is a basis ${\varepsilon }_V$ of $V$ which extends the 
basis ${\varepsilon }_U$ such that the matrix 
of $A$ on $V$ is block diagonal. Let $W$ be the subspace of $V$ with basis 
${\varepsilon }_W = {\varepsilon }_V \setminus {\varepsilon }_U$. 
The matrix of $A|W$ is block diagonal. Therefore $W$ is $A$-invariant 
and $V = U \oplus W$ by construction. Consequently, $A$ is semisimple. 
\hfill $\square $  

\section{The Jordan decomposition of ${\bf A}$}
\label{sec2}

Here we give an algorithm for finding the Jordan decomposition of the linear mapping $A$, 
that is, we find real semisimple and commuting 
nilpotent linear maps $S$ and $N$ whose sum is $A$. The 
algorithm we present uses only the characteristic polynomial 
${\chi }_A$ of $A$ and does {\em not} require that we know 
{\em any} of its factors. Our argument follows that of 
\cite{burgoyne-cushman}. \medskip 

Let $p$ be the square free factorization of ${\chi }_A$. Let $M$ be the 
smallest positive integer such that ${\chi }_A$ divides $p^M$. Then 
$M \le \deg {\chi }_A$. Assume that 
$\deg {\chi }_A \ge 2$, for otherwise $S = A$. Write
\begin{equation}
S = A + \sum^{M-1}_{j=1}r_j(A){p(A)}^j, 
\label{eq-ones2}
\end{equation}
where $r_j$ is a polynomial whose degree 
is less than the degree of $p$. From the fact that ${\chi }_A$ 
divides $p^M$, it follows that ${p(A)}^M =0$.  \medskip 

We want to determine $S$ in the form (\ref{eq-ones2}) so that 
\begin{equation}
p(S) =0.
\label{eq-twos2}
\end{equation}
From claim 1.1 it follows that $S$ is semisimple. \medskip 

We have to find the polynomials $r_j$ in (\ref{eq-ones2}) so 
that equation (\ref{eq-twos2}) holds. We begin by using the Taylor 
expansion of $p$. If (\ref{eq-ones2}) holds, then 
\begin{align}
p(S) & =  p\Big( A + \sum^{M-1}_{j=1}r_j(A){p(A)}^j \Big) \notag \\
& =  p(A) + \sum^{M-1}_{i=1} p^{(i)}(A) \Big( 
\sum^{M-1}_{j=1}r_j(A){p(A)}^j \Big)^i,  \notag \\
& \hspace{.5in}\parbox[t]{3.5in}{where $p^{(i)}$ is $\frac{1}{i!}$ times 
the ${\rm i}^{\rm th}$ derivative of $p$}  \notag \\
& =  p(A) + \sum^{M-1}_{i=1}\sum^{M-1}_{k=1} c_{k,i}\, {p(A)}^k 
p^{(i)}(A). 
\label{eq-twostars2} 
\end{align}
Here $c_{k,i}$ is the coefficient of $z^k$ in 
$(r_1z + \cdots +\, r_{M-1}z^{M-1})^i$. Note that $c_{k,i} =0$ if 
$k>i$. A calculation shows that when $k \le i$ we have 
\begin{equation}
c_{k,i} = \sum_{\stackrel{{\alpha }_1 + \cdots + {\alpha}_{k-1}=i}{{\alpha }_1 + 
2{\alpha }_2 + \cdots + (k-1){\alpha }_{k-1}}= k} \frac{i!}{{\alpha }_1! \cdots 
{\alpha }_{k-1}!} r^{{\alpha }_1}_1 \cdots r^{{\alpha }_{k-1}}_{k-1} .
\label{eq-threestars2}
\end{equation}   
Interchanging the order of summation in (\ref{eq-twostars2}) we get 
\begin{displaymath}
p(S)  =   p(A) + \sum^{M-1}_{i=1}\Big( r_i(A)p^{(1)}(A) + e_i(A) \Big) {p(A)}^i,
\end{displaymath}
where $e_1=0$ and for $i \ge 2$ we have $e_i = \sum^{i}_{j=2} c_{i,j} p^{(j)}$. 
Note that $e_i$ depends on $r_1, \ldots , r_{i-1}$, because of (\ref{eq-threestars2}). \medskip  

Suppose that we can find polynomials $r_i$ and $b_i$ such that 
\begin{equation}
r_ip^{(1)} + e_i = b_i p -b_{i-1},
\label{eq-threes2}
\end{equation}
for every $1 \le i \le M-1$. Here $b_0 =1$. Then 
\begin{displaymath}
\sum^{M-1}_{i=1}\big( r_i(A)p^{(1)}(A) + e_i(A) \big) {p(A)}^i =
\sum^{M-1}_{i=1} \big( b_i(A)p(A) - b_{i-1}(A) \big) {p(A)}^i \, = \, 
-p(A), 
\end{displaymath}
since $p^M(A) =0$ and $b_0 =1$, which implies $p(S) =0$, see (\ref{eq-twostars2}). \medskip 

We now construct polynomials $r_i$ and $b_i$ so that (\ref{eq-threes2}) 
holds. We do this by induction. Since the polynomials $p$ and 
$p^{(1)}$ have no common nonconstant factors, their greatest 
common divisor is the constant polynomial $1$. Therefore by 
the Euclidean algorithm there are polynomials $g$ and $h$ with 
the degree of $h$ being less than the degree of $p$ such that 
\begin{equation}
gp - hp^{(1)} = 1.
\label{eq-fours2}
\end{equation}
\par
Let $r_1 =h$, and $b_1 =g$. Using the fact that $b_0=1$ and $e_1=0$, we see that 
equation (\ref{eq-fours2}) is the same as equation (\ref{eq-threes2}) 
when $i=1$. Let $d_1 = 0$ and $q_0 = q_1 =0$. Now suppose that $n \ge 2$. By induction 
suppose that the polynomials 
$r_1, \ldots \, , r_{n-1}$, $e_1, \ldots , e_{n-1}$, $q_1, \ldots \, , q_{n-1}$ and $b_1, \ldots \, , b_{n-1}$ are known and that $r_i$ and $b_i$ satisfy 
(\ref{eq-threes2}) for every $1 \le i \le n-1$. Using the fact that the polynomials $r_1, \ldots , r_{n-1}$ are known, from formula  
(\ref{eq-threestars2}) we can calculate the polynomial $e_n = \sum^n_{j=2}c_{i,n}\, p^{(j)}$. For  $n \ge 2$ define the polynomial $d_n$ by 
\begin{equation}
d_n = q_{n-1} + h \sum^n_{i=1} g^{n-i} e_i. 
\label{eq-fives2}
\end{equation}
Note that the polynomials $q_{n-1}$, $g = b_1 $, $h= r_1$, and $e_i$ for 
$1 \le i \le n-1$ are already known by the induction hypothesis. Thus the right hand side of (\ref{eq-fives2}) is known and hence so is $d_n$. Now define the polynomials $q_n$ and $r_n$ by dividing $d_n$ by $p$ with remainder, namely
\begin{equation}
d_n = q_n p + r_n .
\label{eq-sixs2}
\end{equation}
Clearly, $q_n$ and $r_n$ are now known. Next for $n\ge 2$ define the 
polynomial $b_n$ by 
\begin{equation}
b_n = -p^{(1)}q_n + g\sum^n_{i=1}g^{n-i}e_i.
\label{eq-sevens2}
\end{equation}
Since the polynomials $p^{(1)}$, $q_n$, $g = b_1$, and $e_i$ for $1 \le i \le n$ are 
known, the polynomial $b_n$ is known. We now show that equation (\ref{eq-threes2}) holds. \medskip 

\noindent \textbf{Proof.} We have already checked that (\ref{eq-threes2}) holds 
when $n=1$. By induction we assumed that it holds for 
every $1 \le i \le n-1$. Using the definition of $b_n$ (\ref{eq-sevens2}) 
and the induction hypothesis we compute 
\begin{align}
b_n p - b_{n-1} & = \Big[ -p^{(1)}pq_n + pg\sum^n_{i=1} g^{n-i}e_i \Big] 
- \Big[ -p^{(1)}q_{n-1} + g\sum^{n-1}_{i=1}g^{n-1-i}e_i \Big] \notag \\
& \hspace{-.5in} =  - p^{(1)}(q_n p - q_{n-1}) + 
pg\sum^n_{i=1}g^{n-i}e_i -\sum^{n-1}_{i=1}g^{n-i}e_i \notag \\
& \hspace{-.5in} =  -p^{(1)}(-r_n + d_n-q_{n-1}) 
+(hp^{(1)}+1)\sum^n_{i=1}g^{n-i}e_i - \sum^{n-1}_{i=1}g^{n-i}e_i, \notag \\
& \hspace{.5in} \mbox{using (\ref{eq-fours2}) and (\ref{eq-sixs2})} \notag  \\
& \hspace{-.5in} =  p^{(1)}r_n -hp^{(1)}\sum^n_{i=1}g^{n-i}e_i +  
hp^{(1)}\sum^n_{i=1}g^{n-i}e_i + \sum^n_{i=1}g^{n-i}e_i - 
\sum^{n-1}_{i=1}g^{n-i}e_i, \notag \\
&\hspace{.5in} \mbox{using (\ref{eq-fives2})} \notag \\
& \hspace{-.5in}=  p^{(1)}r_n +e_n. \tag*{$\square $} 
\end{align}
This completes the construction of the polynomial 
$r_n$ in (\ref{eq-ones2}). Repeating this construction until $n = M-1$ we have determined the semisimple part $S$ of 
$A$. The commuting nilpotent part of A is $N = A-S$. \hfill $\Box $

\section{Uniform normal form}  
\label{sec3}

In this section we give a description of the uniform normal form 
of a linear map $A$ of $V$ into itself. We assume that the Jordan decomposition of $A$ into its commuting semisimple and nilpotent summands $S$ and $N$, respectively, is known. 

\subsection{Nilpotent normal form}
\label{sec3subsec1}

In this subsection we find the Jordan normal form for a nilpotent linear 
transformation $N$. \medskip 

Recall that a linear transformation $N:V \rightarrow V$ is said to be 
{\sl nilpotent of index} $n$ if there is an integer $n \ge 1$ such that 
$N^{n-1} \not = 0$ but $N^n = 0$. Note that the index of nilpotency 
$n$ need not be equal to $\dim V$. Suppose that for some positive integer $\ge 1$ there is a nonzero vector $v$, which lies in $\ker N^{\ell } \setminus \ker N^{{\ell }-1}$. The 
set of vectors $\{ v, Nv, \ldots \, , N^{{\ell }-1}v \} $ is a 
\emph{Jordan chain} of \emph{length} ${\ell }$ with \emph{generating vector} $v$.
The space $V^{\ell }$ spanned by the vectors in a given Jordan chain of 
length ${\ell }$ is a $N$-\emph{cyclic subspace} of $V$. Because 
$N^{\ell }v =0$, the subspace $V^{\ell }$ is $N$-invariant. 
Since $\ker N|V^{\ell } = {\spann }_{\kk } \{ N^{{\ell }-1}v \} $, the mapping 
$N|V^{\ell }$ has exactly one eigenvector corresponding to the 
eigenvalue $0$. \medskip 

\noindent \textbf{Claim 3.1.1} The vectors in a Jordan chain are 
linearly independent. \medskip 

\noindent \textbf{Proof.} Suppose not. Then $0 = \sum^{{\ell }-1}_{i=0} {\alpha }_i\, 
N^iv$, where not every ${\alpha }_i \in \kk $ is zero. Let $i_0$ be the smallest index 
for which ${\alpha }_{i_0} \not = 0$. Then 
\begin{equation}
0  =  {\alpha }_{i_0}\, N^{i_0}v + \cdots \, + {\alpha }_{{\ell }-1} \, 
N^{{\ell }-1}v  .
\label{eq-sec4ss1one}
\end{equation}
Applying $N^{{\ell }-1-i_0}$ to both sides of (\ref{eq-sec4ss1one}) gives 
$0 = {\alpha }_{i_0}N^{{\ell }-1}v$. By hypothesis $v\not \in \ker N^{{\ell }-1}$, 
that is, $N^{{\ell }-1}v \ne 0$. Hence ${\alpha }_{i_0} =0$. This 
contradicts the definition of the index $i_0$. Therefore ${\alpha }_i =0$ for every 
$0 \le i \le \ell -1$. Thus the vectors ${\{ N^iv \}}^{\ell -1}_{i=0}$, which span the 
Jordan chain $V^{\ell }$, are linearly independent. 
\hfill $\square $ \medskip 

With respect to the \emph{standard basis} $\{ N^{{\ell }-1}v, N^{\ell-2}v, 
\ldots \, , Nv, v \}$ of $V^{\ell }$ the matrix of $N|V^{\ell }$ is the ${\ell } \times {\ell }$ matrix 
\begin{displaymath}
\mbox{{\footnotesize $\left( \begin{array}{ccccc} 0 & 1 & 0 & \dots 
& 0 \\
0 & 0 & 1  &\ddots & \vdots \\ 
\vdots & \vdots & \ddots & \ddots & 0 \\
\vdots & \vdots & \vdots & \ddots & 1 \\
0 & 0 & \cdots & \cdots & 0 \\ \end{array} \right) $},}
\end{displaymath}
which is a \emph{Jordan block} of size ${\ell }$. \medskip 

We want to show that $V$ can be decomposed into a 
direct sum of $N$-cyclic subspaces. In fact, we show that there is 
a basis of $V$, whose elements are given by a dark dot $\bullet $ or 
an open dot $\circ $ in the diagram below such that the arrows give 
the action of $N$ on the basis vectors. Such a diagram is called the 
\emph{Young diagram of} $N$.
\begin{displaymath} 
\begin{array}{l}
\phantom{\bullet}\\ 
\phantom{\uparrow}\\ 
\phantom{\bullet}\\ 
\phantom{\uparrow}\\ 
\phantom{\bullet}\\  
\phantom{\vdots}\\ 
\phantom{\bullet}\\ 
\phantom{\uparrow}\\ 
\phantom{circ}\\
\phantom{circ}\\
\phantom{x}\\
\end{array}
\hspace{-1cm} 
\begin{array}{l}
\begin{array}{l}
\phantom{\bullet}\\ 
\phantom{\uparrow}\\ 
\phantom{\bullet}\\ 
\phantom{\uparrow}\\ 
\phantom{\bullet}\\  
\phantom{\vdots}\\ 
\phantom{\bullet}\\ 
\phantom{x}\\
\end{array} \\
\phantom{\uparrow}\\
\phantom{circ}\\
\end{array}
\hspace{-0.5cm}
\begin{array}{l}
\begin{array}{l}
\phantom{\bullet}\\ 
\end{array}\\
\phantom{\uparrow}\\ 
\phantom{\bullet}\\ 
\phantom{\uparrow}\\ 
\phantom{\bullet}\\  
\phantom{\vdots}\\ 
\phantom{\bullet}\\ 
\phantom{\uparrow}\\
\phantom{circ}\\
\end{array}
\hspace{-0.4cm}
\begin{array}{cllllllll} 
\bullet&\bullet&\bullet&\bullet&\circ&\circ&\circ\\
\uparrow&\uparrow&\uparrow&\uparrow\\
\bullet&\bullet&\bullet&\bullet\\ 
\uparrow&\uparrow&\uparrow&\uparrow\\
\bullet&\bullet&\bullet&\circ\\ 
\vdots&\vdots&\vdots\\
\bullet&\bullet&\circ&\\ 
\uparrow&\uparrow\\ 
\circ&\circ 
\end{array}
\end{displaymath}
\par
\noindent
\hspace{1.25in}\parbox[t]{4.5in}{Figure 3.1.1. The Young diagram 
of $N$.} \medskip 

Note that the columns of the Young diagram of $N$ are Jordan chains 
with generating vector given by an open dot. The black dots form a basis 
for the image of $N$, whereas the open dots form a basis for a 
complementary subspace in $V$. The dots on or above the ${\ell }^{\rm th}$ row for a basis for 
$\ker N^{\ell }$ and the black dots in the first row form a basis for $\ker N \cap \mathrm{im}\, N$. Let $r_{\ell }$ be 
the number of dots in the ${\ell }^{\rm th}$ row. Then $r_{\ell } = 
\dim \ker N^{\ell } - \dim \ker N^{\ell -1}$. Thus the Young diagram of $N$ is unique. \medskip

\noindent \textbf{Claim 3.1.2} There is a basis of $V$ that realizes the Young 
diagram of $N$.\medskip 

\noindent \textbf{Proof}. Our proof follows that of Hartl \cite{hartl}. We use induction of the dimension of $V$. Since 
$\dim \ker N > 0$, it follows that $\dim \mathrm{im}\, N < \dim V$. Thus by the induction hypothesis, we may suppose that 
$\mathrm{im}\, N$ has a basis which is the union of $p$ Jordan chains 
$\{ w_i, Nw_i, \ldots \, , N^{m_i}w_i \} $ each of length $m_i$. The \linebreak 
vectors ${\{ N^{m_i}w_i \} }^p_{i=1}$ lie in 
$\mathrm{im}\, N \cap \ker N$ and in fact form a basis of this subspace. 
Since $\ker N$ may be larger than $\mathrm{im}\, N \cap \ker N$, 
choose vectors $\{ y_1, \ldots \, , y_q \} $ where 
$q$ is a nonnegative integer such that $\{ N^{m_1}w_1, \ldots \, , N^{m_p}w_p, $ 
$y_1, \ldots \, , y_q \} $ form a basis of $\ker N$.

Since $w_i \in \mathrm{im}\, N$ there is a vector $v_i$ in $V$ such that 
$w_i = Nv_i$. We assert that the $p$ Jordan chains 
\begin{displaymath}
\{ v_i, Nv_i, \ldots \, , N^{m_i+1}v_i \} = \{ v_i, w_i, Nw_i, 
\ldots \, , N^{m_i}w_i \} 
\end{displaymath}
each of length $m_i +2$ together with the $q$ vectors $\{ y_j \} $, which are  
Jordan chains of length $1$, form a 
basis of $V$. To see that they span $V$, let $v \in V$. Then $Nv 
\in \mathrm{im}\, N$. Using the  
basis of $\mathrm{im}\, N$ given by the induction hypothesis, we may write
\begin{displaymath}
Nv = \sum^p_{i=1}\sum^{m_i}_{\ell =0} {\alpha }_{i\ell } 
N^{\ell } w_i \, = \, 
N \Big( \sum^p_{i=1}\sum^{m_i}_{\ell =0} {\alpha }_{i\ell } 
N^{\ell }v_i \Big).
\end{displaymath}
Consequently,
\begin{displaymath}
v -  \sum^p_{i=1}\sum^{m_i}_{\ell =0} {\alpha }_{i\ell } N^{\ell } 
v_i = \sum^p_{i=1}{\beta }_i N^{m_i+1}v_i + \sum^q_{\ell =1}{\gamma }_{\ell } y_{\ell },
\end{displaymath}
since the vectors
\begin{displaymath}
\{  N^{m_1}w_1, \ldots \, , N^{m_p}w_p, y_1, \ldots \, , y_q \} = 
\{  N^{m_1+1}v_1, \ldots \, , N^{m_p+1}v_p, y_1, \ldots \, , y_q \} 
\end{displaymath}
form a basis of $\ker N$. Linear independence is a consequence of the 
following counting argument. The number of vectors in the Jordan chains is 
\begin{align}
\sum^p_{i=1}(m_i + 2) +q = \sum^p_{i=1}(m_i+1) +(p+q) \, = \, \dim 
\mathrm{im}\, N + \dim \ker N \, = \, \dim V.  \tag*{$\square $} 
\end{align}

We note that finding the generating vectors of the Young diagram of $N$ 
or equivalently the Jordan normal form of $N$, involves solving linear equations with 
coefficients in the field $\kk $ and thus only operations in the field $k$.

\subsection{Some facts about $S$}
\label{sec3subsec3}

We now study the semisimple part $A$. \medskip 

\noindent \textbf{Lemma 3.2.1} $V = \ker S \oplus \mathrm{im}\, S$. 
Moreover the characteristic polynomial ${\chi }_S(\lambda )$ of $S$ can be written as a product of ${\lambda }^n$, where $n = \dim \ker S$ and 
${\chi}_{S|\imm S}$, the characteristic 
polynomial of $S|\imm S$. Note that ${\chi }_{S|\imm S} (0) \ne 0$ \medskip  

\noindent \textbf{Proof.} $\ker S$ is an $S$-invariant subspace of 
$V$. Since $Sv =0$ for every $v \in $ $\ker S$, the characteristic polynomial 
of $S|\ker S$ is ${\lambda }^n$. 

Because $S$ is semisimple, there is an $S$-invariant subspace $Y$ of $V$ 
such that $V = \ker S \oplus Y$. The linear mapping $S|Y:Y \rightarrow Y$ 
is invertible, for if $Sy=0$ for some $y \in Y$, then $S(y+u) =0$ for 
every $u \in \ker S$. Therefore $y+u \in \ker S$, which implies that 
$y \in \ker S \cap Y = \{ 0 \} $, that is, $y=0$. So $S|Y$ is invertible. 
Suppose that $y \in Y$, then $y = S\big( (S|Y)^{-1}y \big) \in \imm S$. 
Thus $Y \subseteq \imm S$. But 
$\dim \imm S = \dim V - \dim \ker S = \dim Y$. So $Y = \imm S$.

Since $\ker S \cap \imm S = \{0 \} $, we see that $\lambda $ does 
not divide the polynomial ${\chi }_{S|\imm S}(\lambda )$. 
Consequently, ${\chi }_{S|\imm S}(0) \ne 0$. Since 
$V = \ker S \oplus \imm S$, where $\ker S$ and $\mathrm{im}\, S$ 
are $S$-invariant subspaces of $V$, we obtain 
\begin{align}
{\chi }_S(\lambda ) & = {\chi }_{\ker S}(\lambda ) \cdot 
{\chi }_{S|\imm S}(\lambda ) = 
{\lambda }^n{\chi }_{S|\imm S}(\lambda ). \tag*{$\square $}
\end{align}    

\noindent \textbf{Lemma 3.2.2} The subspaces $\ker S$ and $\imm S$ 
are $N$-invariant and hence $A$-invariant. \medskip 

\noindent \textbf{Proof.} Suppose that $x \in \imm S$. Then there 
is a vector $v \in V$ such that $x = Sv$. So $Nx = N(Sv) = S(Nv) \in 
\imm S$. In other words, $\imm S$ is an $N$-invariant 
subspace of $V$. Because $\imm S$ is also $S$-invariant and 
$A = S+N$, it follows that $\imm S$ is an $A$-invariant 
subspace of $V$. Suppose that $x \in \ker S$, that is, $Sx =0$. Then 
$S(Nx) = N(Sx) =0$. So $Nx \in \ker S$. Therefore $\ker S$ is an 
$N$-invariant and hence $A$-invariant subspace of $V$. \hfill $\square $ 

\subsection{Decription of uniform normal form}
\label{sec3subsec3}

We now describe the uniform normal form of the linear mapping $A:V \rightarrow V$, 
using both its semisimple and nilpotent parts. \medskip 

Since $A|\ker S = N|\ker S$, we can apply the discussion of \S 3.1 to obtain 
a basis of $\ker S$ which realizes the Young diagram of $N|\ker S$, which say has 
$r$ columns.  For $1 \le {\ell } \le r$ let $F_{q_{\ell }}$ be the space spanned by the generating 
vectors of Jordan chains of $N|\ker S$ in $\ker S$ of length $m_{\ell }$. \medskip 

By lemma 3.2.1 $A|\imm S$ is a linear mapping of $\imm S$ into itself  
with invertible semisimple part $S|\imm S$ and commuting nilpotent part $N|\imm S$. Using 
the discussion of \S 3.1 for every $r+1 \le {\ell } \le p$ let 
$F_{q_{\ell}}$ be the set of generating vectors of the Jordan 
chains of $N|\imm S$ in $\imm S$ of length $m_{\ell }$, which 
occur in the $p-(r+1)$ columns of the Young diagram of $N|\imm S$. \medskip 

Now we prove \medskip 

\noindent \textbf{Claim 3.3.1} For each $1 \le {\ell } \le p$ the space 
$F_{q_{\ell }}$ is $S$-invariant. \medskip 

\noindent \textbf{Proof.} Let $v^{\ell } \in F_{q_{\ell }}$. Then 
$\{ v^{\ell }, \, N v^{\ell }, \ldots , N^{m_{{\ell }}-1}v^{\ell } \} $ is a 
Jordan chain in 
the Young diagram of $N$ of length $m_{\ell }$ with generating vector 
$v^{\ell }$. For each $1 \le {\ell } \le r$ we have $F_{q_{\ell }} \subseteq \ker S$. So trivially $F_{q_{\ell }}$ is $S$-invariant, because $S =0$ on $F_{q_{\ell }}$. Now suppose that 
$r+1 \le {\ell } \le p$. Then $F_{q_{\ell }} \subseteq \imm S $ and 
$S|\imm S$ is invertible  Furthermore, suppose that for some ${\alpha }_j \in \kk $ 
with $0 \le j \le m_{\ell }-1$ we 
have $0 = \sum^{m_{\ell }-1}_{j=0} {\alpha }_j N^j(Sv^{\ell })$. Then 
$0 = S\big( \sum^{m_{\ell }-1}_{j=0}{\alpha }_j N^jv^{\ell } \big)$, because 
$S|\imm S$ and $N|\imm S$ commute. Since $S|\imm S$ is invertible, 
the preceding equality 
implies $0 = \sum^{m_{\ell }-1}_{j=0}{\alpha }_j N^jv^{\ell }$. Consequently, by 
lemma 3.1.1 we obtain ${\alpha }_j =0 $ for every $0 \le j \le m_{\ell }-1$. 
In other words, $\{ Sv^{\ell }, \, N(Sv^{\ell }), \ldots , 
N^{m_{\ell }-1}(Sv^{\ell }) \} $ is 
a Jordan chain of $N|\imm S$ in $\imm S$ of length $m_{\ell }$ with 
generating vector $Sv^{\ell }$. 
So $Sv^{\ell } \in F_{q_{\ell }}$. Thus $F_{q_{\ell }}$ is an $S$-invariant 
subspace of $\imm S$ and hence is an $S$-invariant subspace of $V$, since 
$V = \imm S \oplus \ker S$. \hfill $\square $ \medskip 

An $A$-invariant subspace $U$ of $V$ is \emph{uniform} of \emph{height} $m-1$ if 
$N^{m-1}U \ne \{ 0 \} $ but $N^m U =\{0 \} $ and $\ker N^{m-1}U = NU$. For each $1 \le {\ell } \le r$ let $U^{q_{\ell }}$ be the space spanned by the 
vectors in the Jordan chains of length $m_{\ell }$ in the Young diagram of $N|\ker S$ and 
for $r+1 \le \ell \le p$ let $U^{q_{\ell }}$ be the space spanned by the 
vectors in the Jordan chains of length $m_{\ell }$ in the Young diagram of $N|\imm S$  \medskip

\noindent \textbf{Claim 3.3.2} For each $1 \le {\ell } \le p$ the subspace  
$U^{q_{\ell }}$ is uniform of height $m_{\ell }-1$. \medskip

\noindent \textbf{Proof.} By definition $U^{q_{\ell }} = F_{q_{\ell }} \oplus NF_{q_{\ell }} 
\oplus \cdots \oplus N^{m_{\ell }-1}F_{q_{\ell }}$. Since 
$N^{m_{\ell }}F_{q_{\ell }} =\{ 0 \} $ but $N^{m_{\ell }-}F_{q_{\ell }} \ne \{ 0 \} $, 
the subspace $U^{q_{\ell }}$ is $A$-invariant and of 
the height $m_{\ell }-1$. To show that $U^{q_{\ell }}$ is uniform we need only show that  
$\ker N^{m_{\ell }-1} \cap U^{q_{\ell }} \subseteq NU^{q_{\ell }}$ since the inclusion 
of $NU^{q_{\ell }}$ in $\ker N^{m_{\ell }-1}$ follows from the fact that 
$N^{m_{\ell }}F_{q_{\ell }} =0$. Suppose that $u \in 
\ker N^{m_{\ell }-1} \cap U^{q_{\ell }}$, 
then for every $0 \le i \le m_{\ell }-1$ there are unique vectors 
$f_i \in F_{q_{\ell }}$ such that $u = f_0 + Nf_1 + \cdots 
+ N^{m_{\ell }-1}f_{m_{\ell }-1}$. Since $u \in \ker N^{m_{\ell }-1}$ we get 
$0 = N^{m_{\ell }-1}u = N^{m_{\ell }-1}f_0$. If $f_0 \ne 0$, 
then the preceding equality contradicts the fact that $f_0$ is a generating vector of 
a Jordan chain of $N$ of length $m_{\ell }$. Therefore $f_0 =0$, 
which means that $u = N(f_1 + \cdots + N^{m_{\ell }-2}f_{m_{\ell }-1}) \in NU^{q_{\ell }}$. This shows that $\ker N^{m_{\ell }-1} \cap U^{q_{\ell }} \subseteq NU^{q_{\ell }}$. Hence 
$\ker N^{m_{\ell }-1} \cap U^{q_{\ell }} = NU^{q_{\ell }}$, that is, the 
subspace $U^{q_{\ell }}$ is uniform of height $m_{\ell }-1$. \hfill $\square $ \medskip     

Now we give an explicit description of the uniform normal form of the linear mapping $A$. For each $1 \le {\ell } \le p$ let ${\chi }_{S|F_{q_{\ell }}}$ be the characteristic polynomial of $S$ on $F_{q_{\ell }}$. From the fact that every summand in 
$U^{q_{\ell }} = F_{q_{\ell }} \oplus NF_{q_{\ell }} \oplus \cdots \oplus N^{m_{\ell }-1}F_{q_{\ell }}$ 
is $S$-invariant, it follows that the 
characteristic polynomial ${\chi }_{S|U^{q_{\ell }}}$ of $S$ on $U^{q_{\ell }}$ is 
${\chi }^{m_{\ell }}_{S|F_{q_{\ell }}}$. Since $V = \sum^p_{\ell =1}\oplus U^{q_{\ell }}$, we 
obtain ${\chi }_S = \prod^p_{\ell =1}{\chi }^{m_{\ell }}_{S|F_{q_{\ell }}}$. 
Choose a basis ${\{ u^{\ell }_j \} }^{q_{\ell }}_{j=1}$ 
of $F_{q_{\ell }}$ so that the matrix of $S|F_{q^{\ell }}$ is the $q_{\ell } \times q_{\ell }$ 
companion matrix $C_{q_{\ell }}$ (\ref{eq-twostars1}) associated to the characteristic polynomial ${\chi }_{S|F_{q_{\ell }}}$. When $1 \le {\ell } \le r$ the companion matrix 
$C_{q_{\ell }}$ is $0$ since $S|F_{q_{\ell }} =0$. With respect to the basis 
${\{ u^{\ell }_j, \, Nu^{\ell }_j, 
\ldots , N^{m_{\ell }-1}u^{\ell }_j \} }^{q_{\ell }}_{j=1} $ of $U^{q_{\ell }}$ 
the matrix of $A|U^{q_{\ell }}$ is the $m_{\ell }q_{\ell } \times m_{\ell }q_{\ell }$ matrix 
\begin{displaymath}
D_{m_{\ell }q_{\ell }} = \mbox{\footnotesize $\left( \begin{array}{cccccl}
C_{q_{\ell }}  &   0      &  0      & \cdots  &  \cdots  & 0       \\
I        & C_{q_{\ell }}  &  0      & \cdots  &  \vdots  & 0       \\ 
0        &  I       & \ddots  &         &  \vdots  & \vdots  \\
\vdots   &          &  \ddots &  \ddots  & \vdots  & \vdots  \\ 
0        &  \cdots  &  0      &   I     & C_{q_{\ell }}  & 0       \\ 
0        & \cdots   &  \cdots &   0     & I     & C_{q_{\ell }}  
\end{array} \right) .$}
\end{displaymath}
Since $V = \sum^p_{{\ell }=1} \oplus U^{q_{\ell }}$, the matrix of $A$ is $\mathrm{diag}\, (D_{m_1q_1}, \ldots , D_{m_pq_p})$ with respect to the basis 
${\{ u^{\ell }_j, \, Nu^{\ell }_j, \ldots , N^{m_{\ell }-1}u^{\ell }_j \} }^{(q_{\ell }, p)}_{(j, \ell ) = (1,1)}$.  
We call preceding matrix the \emph{uniform normal form} for the linear map $A$ of $V$ 
into itself. We note that this normal form can be computed using only operations in the field 
$\kk $ of characteristic $0$.  \medskip 

Using the uniform normal form of $A$ we obtain a factorization of its characteristic polynomial 
${\chi }_A$ over the field $\kk $. \medskip 

\noindent \textbf{Corollary 3.3.3} ${\chi }_A (\lambda )=  
\prod^p_{\ell =1}{\chi }^{m_{\ell }}_{S|F_{q_{\ell }}}(\lambda ) = {\lambda }^n \, 
\prod^p_{\ell = r+1}{\chi }^{m_{\ell }}_{S|F_{q_{\ell }}}(\lambda )$, where 
$n = \sum^r_{\ell =1} m_{\ell } = \dim \ker S$.

\end{document}